\documentclass[a4paper,11pt,reqno,twoside]{article}
\usepackage{amssymb}
\usepackage[frenchb,english]{babel}
\usepackage{amscd}
\usepackage{amsmath,amsthm,amsfonts,amssymb,graphicx}
\usepackage[colorlinks,linkcolor=blue,anchorcolor=blue,citecolor=green]{hyperref}
\usepackage{mathrsfs}
\usepackage{fancyhdr}
\usepackage{multicol}
\usepackage{abstract}
\usepackage{mathrsfs}

\begin{document}
\theoremstyle{plain}
\newtheorem{Definition}{Definition}[section]
\newtheorem{Proposition}{Proposition}[section]
\newtheorem{Property}{Property}[section]
\newtheorem{Theorem}{Theorem}[section]
\newtheorem{Lemma}{Lemma}[section]
\newtheorem{Corollary}[Theorem]{Corollary}
\newtheorem{Remark}{Remark}[section]

\noindent  {\LARGE
K\"{a}hler submanifolds of the symmetrized polydisc }\\\\
\noindent\text{Guicong Su, \; Yanyan Tang \; \& \; Zhenhan Tu$^{*}$ }\\

\noindent\small {School of Mathematics and Statistics, Wuhan
University, Wuhan, Hubei 430072, P.R. China}

\noindent\text{Email: suguicong@whu.edu.cn (G. Su), yanyantang@whu.edu.cn (Y. Tang), zhhtu.math@whu.edu.cn (Z. Tu)}

\renewcommand{\thefootnote}{{}}
\footnote{\hskip -16pt {$^{*}$Corresponding author. \\}}\\

\normalsize \noindent\textbf{Abstract}\quad   This paper proves the non-existence of common K\"{a}hler submanifolds of the complex Euclidean space and the symmetrized polydisc endowed with their canonical metrics.

\vskip 10pt

\vskip 10pt


\pagenumbering{arabic}
\renewcommand{\theequation}
{\arabic{section}.\arabic{equation}}

\setcounter{section}{0}
\setcounter{equation}{0}
\section{Introduction}

The study of holomorphic isometric embeddings between complex manifolds is a classical topic in complex geometry. For K\"{a}hler manifolds, the early work was carried out by Calabi who, in his seminal work \cite{Cal} in 1953, proved that two complex space forms with different curvature signs cannot be locally isometrically embedded into another one with respect to the canonical K\"{a}hler metrics. In the proof of his result, the notion $``diastasis"$ for analytic  K\"{a}hler metrics plays an essential role.

Along these lines of research, Umehara \cite{Ume2} proposed an interesting question whether two K\"{a}hler manifolds have in common a non-trivial K\"{a}hler submanifold with the induced metrics, and showed that K\"{a}hler submanifolds of complex space form of different types are essentially different from each other. Inspired by the work of Umehara, Di Scala and Loi \cite {DL} introduced the concept of $``relatives"$ between two K\"{a}hler manifolds (i.e., they are said to be relatives if they share a common K\"{a}hler submanifold, otherwise, we say that they are not relatives) in 2010, and they proved that a bounded domain with its Bergman metric and a projective K\"{a}hler manifold with the restriction of the Fubini-Study metric are not relatives. For related problems, see  Cheng, Di Scala and Yuan \cite{CDY},  Di Scala and Loi \cite{DL0},   Mossa \cite{Mos}  and Zedda \cite{Z}.

Notice that any irreducible Hermitian symmetric space of compact type can be holomorphically isometrically embedded into a complex project space by the classical Nakagawa-Takagi embedding.
Therefore, associating Umehara's main result in \cite{Ume2} with the Nakagawa-Takagi embedding, it is easy to get that complex Euclidean spaces and Hermitian symmetric spaces of compact types are not relatives. What is more, Huang and Yuan \cite{HY1} solved the problem of the non-relativity between complex Euclidean spaces and the bounded symmetric domains by using the properties of Nash algebraic function.

More recently, Cheng and Niu \cite{Cheng} discussed the non-relativity in the nonhomogeneous setting.
Cartan-Hartogs domains, introduced by Yin and Roos, are defined as the Hartogs type domain constructed over the bounded symmetric domains. They are natural generalizations of bounded symmetric domains and ellipsoids, but in general they are not homogeneous (e.g., see Feng-Tu \cite{FT},  Yin \cite{Y}). Cheng and Niu \cite{Cheng} studied the non-relativity of a complex Euclidean space and a Cartan-Hartogs domain with canonical metrics.

We define the symmetrized polydisc $\mathbb{G}_n$ as follows. Let $\mathbb{D}$ be the unit disc in the complex plane $\mathbb{C}$, $\lambda=(\lambda_1, \ldots, \lambda_n)\in \mathbb{C}^{n}$, and let $\pi_n=(\pi_{n,1}, \ldots, \pi_{n,n}): \mathbb{C}^{n}\longrightarrow \mathbb{C}^{n}$ be the symmetrization map defined by
\begin{equation} \label{10}
\pi_{n,k}(\lambda)=\sum_{1\leq j_1 < \cdots < j_k \leq n}\lambda_{j_1} \ldots \lambda_{j_k},~~1\leq k \leq n.
\end{equation}
The image $\mathbb{G}_n:=\pi_{n}(\mathbb{D}^{n})$ is known as the symmetrized polydisc. In particular,  $G_1=\mathbb{D},$   and $G_2$ is the so-called symmetrized bidisc.
For the general reference of symmetrized polydisc, see Edigarian-Zwonek \cite{Ber} and the Chapter 7 of Jarnicki-Pflug's book \cite{JP1}.

The symmetrized polydisc $\mathbb{G}_n$ $(n \geq 2)$  is a bounded inhomogeneous pseudoconvex domain without smooth boundary, and especially it hasn't any strongly pseudoconvex boundary point. It is important because the symmetrized bidisc is the first known example of a bounded pseudoconvex domain for which
the Lempert function, the Kobayashi distance and the Carath\'eodory distance coincide, but which cannot be exhausted by domains biholomorphic to convex ones (see Costara \cite{CK}).
The symmetrized polydisc has been studied by many authors, especially in 2-dimensional case, e.g., Agler-
Lykova-Young \cite{ALY}, Agler-Young \cite{AY, starlike}, Frosini-Vlacci \cite{FV2}, Tu-Zhang \cite{TU} and
Trybula \cite{Inv}.

The restriction map $\pi_{n}|_{\mathbb{D}^{n}}: \mathbb{D}^{n}\rightarrow \mathbb{G}^{n}$ is a
proper holomorphic map (see Rudin \cite{Rudin}). Thus, the symmetrized polydisc $\mathbb{G}_n$ is a proper image of the bounded symmetric domain $\mathbb{D}^{n}$.
The purpose of this paper is to prove that the non-existence of common K\"{a}hler submanifolds of the complex Euclidean space and the symmetrized polydisc endowed with their canonical metrics.
Denote Euclidean metric on the complex Euclidean space $\mathbb{C}^m$ and Bergman metric on the symmetrized polydisc $\mathbb{G}_n$ by $\omega_{\mathbb{C}^m}$ and $\omega_{\mathbb{G}_n}$, respectively. Let $\omega_D$ be a K\"{a}hler metric (not necessarily complete) on a domain $D \subseteq \mathbb{C}$ (assume without loss of generality 0 is in $D$).
In this paper, we show that there do not simultaneously exist holomorphic isometric
immersions  $F: (D, \omega_D)\rightarrow (\mathbb{C}^m,  \omega_{\mathbb{C}^m})$ and $G: (D, \omega_D) \rightarrow (\mathbb{G}_n, \omega_{\mathbb{G}_n})$ as follows.

\begin{Theorem}\label{Th1.1}
Let $D$ be a domain in $\mathbb{C}$. Assume that $F=(f_1, \ldots, f_m): D\rightarrow \mathbb{C}^m$ and $G=(g_1, \ldots, g_n): D\rightarrow \mathbb{G}_n$ are holomorphic mappings such that
\begin{equation}\label{1.1}
F^{*}\omega_{\mathbb{C}^m}=G^{*}\omega_{\mathbb{G}_n} \quad on~D.
\end{equation}
Then $F$ must be a constant map.
\end{Theorem}

As an immediate consequence, we have the following:

\begin{Corollary}
There does not exist a K\"{a}hler manifold $(X, \omega_X)$ that can be holomorphic isometrically
embedded into the complex Euclidean space $(\mathbb{C}^m,  \omega_{\mathbb{C}^m})$ and also into the symmetrized polydisc $(\mathbb{G}_n, \omega_{\mathbb{G}_n})$.
\end{Corollary}

Our proof uses the idea developed in the work of Huang and Yuan \cite{HY1}, but due to the fact that the nonhomogeneous of $\mathbb{G}_n$  $(n \geq 2)$, we cannot assume that $G(0)=0$ just as their proof without loss of generality. On the other hand, although that some Cartan-Hartogs domains, researched by Cheng and Niu \cite{Cheng} with the same problem, are also inhomogeneous, they \cite{Cheng} stressed the extra condition $G(0)=0$ there.
We make no such assumption about ${G}$ in the paper. And then, this causes that  $K_{\mathbb{G}_n}\big(G(z), {G}(0)\big)$ may not be a constant on $D$ in the process of our proof, where $K_{\mathbb{G}_n}(\cdot, \cdot)$ is the Bergman kernel of $\mathbb{G}_n$. The key ideas in this paper is to verify that the Bergman kernel  $K_{\mathbb{G}_n}(\cdot, \cdot)$ of $\mathbb{G}_n$ can be described as a rational form and
$\{f_i(z)\}_{i=1}^m$ can be written as holomorphic rational functions in $g_1(z), \ldots, g_n(z)$ on $D$.

\section{Preliminaries}

In this section we review several basic facts about the  symmetrized polydisc $\mathbb{G}_n$, the fundamental theorem of symmetric polynomials, and Nash-algebraic function, which will be used in the subsequent section.

\setcounter{equation}{0}
\renewcommand{\theequation}{2.\arabic{equation}}

For the proper holomorphic map $\pi_{n}: \mathbb{D}^{n}\rightarrow \mathbb{G}^{n}\;(n\geq 2)$,
we have
$$\det\pi'_n(\lambda)= \prod_{1\leq i<j\leq n}(\lambda_i- \lambda_j).$$
Define
$$\Sigma_n:=\{\lambda \in \mathbb{D}^{n}: \det \pi_{n}'(\lambda)=0\}.$$
Then
\begin{equation}\label{200}
{\pi_{n}}^{-1}(\pi_{n}(\lambda)) = \{ (\lambda_{\sigma(1)},\cdots, \lambda_{\sigma(n)}): \sigma\; \mbox{is a permutation of}\; \{1,\cdots,n\}\},\;\;  \lambda\in \mathbb{D}^{n}\setminus \Sigma_n.
\end{equation}
The proper holomorphic map $\pi_{n}: \mathbb{D}^{n}\rightarrow \mathbb{G}^{n}  $ is also a branched covering.
For a holomorphic mapping $G: D \rightarrow \mathbb{G}_n$ of a simply connected domain $D \subseteq \mathbb{C}$ with $G(D)\bigcap  \pi_{n}(\Sigma_n)\not= \emptyset$ (specially, $G(D)\subset\pi_{n}(\Sigma_n)$), generally speaking, it may be impossible to find a holomorphic mapping $T: D \rightarrow \mathbb{D}^n$ such that $G\equiv\pi_{n}(T)$ on $D$.

Let $K_{\Omega}(\cdot, \cdot)$ denote the Bergman kernel of the domain $\Omega \subseteq \mathbb{C}^n$. From the formula for the Bergman kernel of the polydisc $\mathbb{D}^{n}$ and from the formula
for the behavior of the Bergman kernel under proper holomorphic mappings (see Bell \cite{Bell}), by (\ref{200}), Edigarian-Zwonek \cite{Ber} obtained the following result.

\begin{Lemma} $(\mathrm{See}\,\emph{\cite{Ber}})$\label{lem2.1}
\begin{equation}\label{2.1}
K_{\mathbb{G}_{n}}\big(\pi_n(\lambda), {\pi_n(\mu)}\big)=\frac{\det\big[\frac{1}{(1-\lambda_{j}\bar{\mu}_{k})^2}\big]_{1 \leq j,k \leq n}}{\pi^{n}\det \pi_{n}'(\lambda)\det\overline{ \pi_{n}'(\mu)}},~~\lambda, \mu \in \mathbb{D}^{n}\backslash \Sigma_n.
\end{equation}
\end{Lemma}

Observe that although the right-hand side of (\ref{2.1}) is not formally defined on the whole $\mathbb{G}_n\times \mathbb{G}_n$, it extends smoothly on this set.
In the case $n=2$, the elementary calculation shows the following.

\begin{Lemma} $(\mathrm{See}\,\emph{\cite{Ber}})$\label{lem2.2}
\begin{eqnarray*}
K_{\mathbb{G}_{2}}\big(\pi_2(\lambda), {\pi_2(\mu)}\big)=\frac{2-(\lambda_1+\lambda_2)(\bar{\mu}_1+\bar{\mu}_2)+
2\lambda_1\lambda_2\bar{\mu}_1\bar{\mu}_2}{\pi^2[(1-\lambda_1\bar{\mu}_1)
(1-\lambda_1\bar{\mu}_2)(1-\lambda_2\bar{\mu}_1)(1-\lambda_2\bar{\mu}_2)]^2}.
\end{eqnarray*}
\end{Lemma}

Using Lemma \ref{lem2.2}, we easily obtain
\begin{equation*}
K_{\mathbb{G}_{2}}\big((s_1, p_1), ({s}_{2}, {p}_{2})\big)
=\frac{2-s_1\bar{s}_{2}+2p_1\bar{p}_{2}}{\pi^2[1-s_1\bar{s}_{2}+(s_1^2-2p_1)\bar{p}_{2}
-p_1s_1\bar{s}_{2}\bar{p}_{2}+p_1\bar{s}_{2}^2+p_1^2\bar{p}_{2}^2]^2},
\end{equation*}
where $(s_i, p_i)\in \mathbb{G}_{2}, i=1, 2$.

Then the above equation gives an explicit formula  for $K_{\mathbb{G}_2}(\cdot, \cdot)$, which is independent of the symmetrization map $\pi_2$. Moreover, $K_{\mathbb{G}_2}(\cdot, \cdot)$ is a rational function on $\mathbb{G}_2\times \mathbb{G}_2$. However, it seems difficult to write an explicit formula for $K_{\mathbb{G}_n}(\cdot, \cdot)$ $(n\geq3)$ (e.g., see Remark 12 in Edigarian-Zwonek \cite{Ber}).

In order to obtain a more handy form for $K_{\mathbb{G}_{n}}(\cdot, \cdot)$,  we need the fundamental theorem of symmetric polynomials as follows.
\begin{Lemma} $(\mathrm{See}\,\emph{\cite{jibendingli2}})$\label{lem2.3}
Any symmetric polynomial in $n$ variables $x_1,\ldots, x_n$ is representable in a unique way as a polynomial in the $n$ elementary symmetric polynomials $\sigma_1, \ldots, \sigma_n$, where $\sigma_j$ is the $j$th elementary symmetrized polynomial, i.e.,
\begin{equation*}
\sigma_j=\sum_{1\leq i_1<\cdots <i_j\leq n}x_{i_1} \cdots x_{i_j}.
\end{equation*}
\end{Lemma}

Next let us recall some properties of Nash-algebraic function. A holomorphic function $F$ over $U \subseteq \mathbb{C}^k$ is called a holomorphic Nash-algebraic function if there is a non-zero holomorphic polynomial $P(z, X)$ in $X$ with coefficients in polynomials of $z$ such that $P(z, F(z))\equiv 0$ on $U$.
Furthermore, one can assume that $P(z, X)$ is an irreducible polynomial
\begin{equation*}
P(z, X)=a_d(z)X^d+a_{d-1}(z)X^{d-1}+\cdots+a_0(z)
\end{equation*}
where $a_i$ $(i=0, \ldots, d)$ are holomorphic polynomials in $z$ having no common factors and $a_d\not\equiv0$. $P(z, X)$ is said to be the annihilating function of $F(z)$.
 Huang and Yuan \cite{HY1} obtained the following result.

\begin{Lemma}  $(\mathrm{See}\,\emph{\cite{HY1}})$\label{lem2.4}
Let $U \subseteq \mathbb{C}^k$ be a connected open set, and $\xi=(\xi_1, \ldots, \xi_k) \in U$. Let $H_1(\xi), \ldots, H_l(\xi)$ and $H(\xi)$ be holomorphic Nash-algebraic functions on $U$. Assume that
\begin{equation*}
\exp{H(\xi)}=\prod_{i=1}^{l}(H_i(\xi))^{\mu_i}
\end{equation*}
for certain real numbers $\mu_1, \ldots, \mu_l$. Then $H(\xi)$ is constant on $U$.
\end{Lemma}

\section{Proof of the main results}

\setcounter{equation}{0}
\renewcommand{\theequation}{3.\arabic{equation}}

Assume, to reach a contradiction, that $F: D\rightarrow \mathbb{C}^n$ is not constant. Assume without loss of generality that $D$ is simply connected,  0 is in $D$, and $F(0)=0$.

Using the condition (\ref{1.1}), we obtain
\begin{equation}\label{3.1}
\partial\bar{\partial}\Big(\sum_{i=1}^{m}|f_i(z)|^2-\log K_{\mathbb{G}_n}\big(G(z), {G(z)}\big)\Big)=0,\ \ \  z \in D.
\end{equation}

This means that the term  $\sum_{i=1}^{m}|f_i(z)|^2-\log K_{\mathbb{G}_n}\big(G(z), {G(z)}\big)$ in the equation (\ref{3.1}) is a harmonic function of $z$ on the domain $D$ in $\mathbb{C}$. So there exists a holomorphic function $\varphi(z)$ on $D$ such that
\begin{equation}\label{3.0}
\sum_{i=1}^{m}|f_i(z)|^2-\log K_{\mathbb{G}_n}\big(G(z), {G(z)}\big)=\varphi(z)+\overline\varphi(z),\ \ \  z \in D.
\end{equation}
By $F(0)=0$, we have $2\text{Re}\varphi(0)=-\log K_{\mathbb{G}_n}\big(G(0), {G(0)}\big)$. Then by complexifying (\ref{3.0}), we get
\begin{equation}\label{3.2}
\sum_{i=1}^{m}f_i(z)\bar{f_i}(w)-\log K_{\mathbb{G}_n}\big(G(z), {G}(w)\big)=\varphi(z)+\overline\varphi(w),
\end{equation}
where $(z, w)\in D \times D$. Hence let $w=0$, we obtain
\begin{equation}\label{har}
\varphi(z)=-\log K_{\mathbb{G}_n}\big(G(z), {G}(0)\big)-\overline\varphi(0).
\end{equation}

Next we divide into three steps to prove Theorem \ref{Th1.1}.

\textbf{Step 1.} We claim that the Bergman kernel $K_{\mathbb{G}_n}(\xi, {\eta})$ is a rational function in $\xi$ and $\bar{\eta}$.

Denote $\xi=\pi_n(\lambda),\eta=\pi_n(\mu)$. If $\lambda, \mu \in \mathbb{D}^n\backslash\Sigma_n$, then, by the calculating of the numerator in the formula (\ref{2.1}), we get
\begin{equation}\label{ker}
\begin{aligned}
K_{\mathbb{G}_{n}}(\xi, {\eta})&=
K_{\mathbb{G}_{n}}\big(\pi_n(\lambda), {\pi_n(\mu)}\big)\\
&=\frac{\det\big[\frac{1}{(1-\lambda_{j}\bar{\mu}_{k})^2}\big]_{1 \leq j,k \leq n}}{\pi^{n}\det \pi_{n}'(\lambda)\det\overline{ \pi_{n}'(\mu)}}\\
&=\frac{1}{\prod_{i,j=1}^{n}
(1-\lambda_i\bar{\mu}_j)^2}\frac{P_1(\lambda, \bar{\mu})}{
\det \pi_{n}'(\lambda)\det\overline{ \pi_{n}'(\mu)}},
\end{aligned}
\end{equation}
where $P_1(\lambda, \bar{\mu})$ is a polynomial in $\lambda$ and $\bar{\mu}$,  and $\det \pi_{n}'(\lambda) =\prod\nolimits_{1\leq j< k \leq n}(\lambda_j-\lambda_k).$

Let $P_2(\lambda, \bar{\mu}):=\prod_{i,j=1}^{n}
(1-\lambda_i\bar{\mu}_j)^2$. Notice that $P_2(\lambda, \bar{\mu})$ is zero-free on $\mathbb{D}^n\times\mathbb{D}^n$ and
\begin{eqnarray}\label{3.4}
P_2(\lambda, \bar{\mu})&=&\prod_{j=1}^{n}
\big[(1-\lambda_1\bar{\mu}_j)\cdots (1-\lambda_n\bar{\mu}_j)\big]^2 \nonumber \\
&=&\prod_{j=1}^{n}
(1-\xi_1\bar{\mu}_j+\xi_2 \bar{\mu}_{j}^2+\cdots+(-1)^n\xi_n\bar{\mu}_{j}^n)^2\\
&:=&\widetilde{P}_{\xi}(\bar\mu),\nonumber
\end{eqnarray}
where  $\xi=(\xi_1,\ldots, \xi_n):=\pi_n(\lambda).$
Since $\widetilde{P}_{\xi}(\bar\mu)$ is a symmetric polynomial in $\bar{\mu}$, by applying Lemma \ref{lem2.3} to equation (\ref{3.4}), we have
\begin{equation}\label{3.5}
\widetilde{P}_{\xi}(\bar\mu)=\sum_{\alpha}h_{\alpha}(\xi)\bar{\eta}^{\alpha}\;(=\sum_{\alpha}h_{\alpha}(\xi)\overline{\pi_n(\mu)}^{\alpha}),
\end{equation}
where the sum over $\alpha$ with finite terms. Together with (\ref{3.4}) and (\ref{3.5}), we conclude that $h_{\alpha}(\xi)$ is a holomorphic polynomial in $\xi$ by the uniqueness of the power series. Hence,
$$P_2(\lambda, \bar{\mu})=H_2(\xi, \bar{\eta})$$
for a polynomial $H_2(\xi, \bar{\eta})$ in $\xi$ and $\bar{\eta}$. Notice that $H_2(\xi, \bar{\eta})$ is holomorphic in $\xi$ and anti-holomorphic in $\eta$, and what's more, it is zero-free on $\mathbb{G}_n\times\mathbb{G}_n$.

Now we consider the remaining terms in (\ref{ker}).  Let
$$\widetilde{P}_1(\lambda, \bar{\mu}):=\frac{P_1(\lambda, \bar{\mu})}{\det \pi_{n}'(\lambda)\det\overline{ \pi_{n}'(\mu)}}.$$
Since the  Bergman kernel $K_{\mathbb{G}_{n}}\big(\xi, \eta\big)$ extends smoothly on the whole $\mathbb{G}_n\times\mathbb{G}_n$, this means that
$\widetilde{P}_1(\lambda, \bar{\mu})$
is a rational function in $\lambda$ and $\bar{\mu}$, and  smoothly on $\mathbb{D}_n\times\mathbb{D}_n.$
Note
\begin{equation}\label{3.x}
\det \pi_{n}'(\lambda) \det\overline{\pi_{n}'(\mu)} =\prod\nolimits_{1\leq j< k \leq n}(\lambda_j-\lambda_k) \overline{\prod\nolimits_{1\leq r< s \leq n}(\mu_r-\mu_s)}.
\end{equation}
We claim that $\widetilde{P}_1(\lambda, \bar{\mu})$ must be a polynomial in $\lambda$ and $\bar\mu$.
In fact, let $\widetilde{P}_1(\lambda, \bar{\mu}):=\frac{\phi_1(\lambda,\bar\mu)}{\phi_2(\lambda,\bar\mu)}$, where $\phi_1$ and $\phi_2$ are polynomial on $\mathbb{C}^n\times\mathbb{C}^n$ with
$\dim_{\mathbb{C}}\phi_1^{-1}(0)\cap  \phi_2^{-1}(0)\leq 2n-2.$    Since $\phi_2(\lambda,\bar\mu)$ is a factor of the polynomial $\det \pi_{n}'(\lambda)\det\overline{ \pi_{n}'(\mu)}$, we have
 $\phi_2$ must be a constant. Otherwise by (\ref{3.x}), $\phi_2(\lambda,\bar\mu)$ has complex $(2n-1)-$dimensional zero set in  $\mathbb{D}_n\times\mathbb{D}_n$, a contradiction with  $\widetilde{P}_1(\lambda, \bar{\mu})$
smoothly on $\mathbb{D}_n\times\mathbb{D}_n.$

Thus, we have
\begin{equation*}
\widetilde{P}_1(\lambda, \bar{\mu})=K_{\mathbb{G}_{n}}(\xi, {\eta})H_2(\xi, \bar{\eta})
\end{equation*}
for $(\xi, \eta)=(\pi_n(\lambda), \pi_n(\mu))\in\mathbb{G}_n\times\mathbb{G}_n$, and $(\lambda, \mu) \in \mathbb{D}^n\times\mathbb{D}^n$. This implies that there exists a real analytic function $H_1(\xi, \bar{\eta})$ on $\mathbb{G}_n\times\mathbb{G}_n$ , such that $\widetilde{P}_{1}(\lambda, \bar{\mu})=H_1(\xi, \bar{\eta})$. Since $\widetilde{P}_1(\lambda, \bar{\mu})$ is a polynomial in $\lambda$ and $\bar\mu$, we have that $H_1(\xi, \bar{\eta})$ is also a polynomial in $\xi$ and $\bar{\eta}$.

Therefore, the Bergman kernel $K_{\mathbb{G}_n}(\cdot, \cdot)$ can be written as follows:
\begin{equation*}
K_{\mathbb{G}_{n}}(\xi, {\eta})=\frac{H_1(\xi, \bar{\eta})}{H_2(\xi, \bar{\eta})}, ~~(\xi, \eta)\in \mathbb{G}_n\times\mathbb{G}_n,
\end{equation*}
i.e., $K_{\mathbb{G}_n}(\cdot, \cdot)$ is a rational function.

\textbf{Step 2.} We claim that $\{f_i(z)\}_{i=1}^m$ can be written as holomorphic rational functions in $g_1(z), \ldots, g_n(z)$, shrinking $D$ towards the origin if needed.

Since $K_{\mathbb{G}_{n}}(\cdot, \cdot)$ is a real analytic function with $K_{\mathbb{G}_n}(\xi, {\xi})\neq0$ on $\mathbb{G}_{n}\times \mathbb{G}_{n}$, we assume that there exist two neighborhoods $U_1$ and $U_2$ of the origin such that $K_{\mathbb{G}_n}\big(G(z), {G}(w)\big)\neq0$ for all $(z,w)\in U_1\times U_2$, and it also holds for $H_1\big(G(z), \bar{G}(w)\big)$.

Now we fix $z$ near 0,  and differentiate formula (\ref{3.2}) with respect to $\bar w$ near 0. Then, by writing
\begin{equation}\label{3.y}
 D^{\delta}\big(\bar{F}(w)\big) = \big(\frac{\partial^\delta}{\partial \bar w^{\delta}}\bar{f}_{1}(w), \ldots, \frac{\partial^\delta}{\partial \bar w^{\delta}}\bar{f}_n(w)\big),~~\forall \delta \in \mathbb{N},
\end{equation}
we obtain
\begin{eqnarray*}
F(z) \cdot D^{1}\big(\bar{F}(w)\big) = \frac{1}{K_{\mathbb{G}_n}\big(G(z), {G}(w)\big)}\bigg\{\sum_{k=1}^{n}\frac{\partial K_{\mathbb{G}_n}}{\partial \bar{\eta}_k}\big(G(z), {G}(w)\big)\frac{\partial \bar{g}_k}{\partial \bar w}(w)\bigg\}+\overline\varphi'(w).
\end{eqnarray*}
If we set $w=0$ in this formula, simple calculation shows
\begin{eqnarray}\label{3.7}
F(z) \cdot D^{1}\big(\bar{F}(0)\big)& = &\frac{H_2\big(G(z), \bar{G}(0)\big)\sum_{k=1}^{n}\frac{\partial H_1}{\partial \bar{\eta}_k}\big(G(z), \bar{G}(0)\big)\frac{\partial \bar{g}_k}{\partial \bar w}(0)}{K_{\mathbb{G}_n}\big(G(z), {G}(0)\big)\big(H_2(G(z), \bar{G}(0))\big)^2} \nonumber\\
&-&\frac{H_1\big(G(z), \bar{G}(0)\big)\sum_{k=1}^{n}\frac{\partial H_2}{\partial \bar{\eta}_k}\big(G(z), \bar{G}(0)\big)\frac{\partial \bar{g}_k}{\partial \bar w}(0)}{K_{\mathbb{G}_n}\big(G(z), {G}(0)\big)\big(H_2(G(z), \bar{G}(0))\big)^2} \nonumber+\overline\varphi'(0)\\
&=&\frac{\sum_{k=1}^{n}\frac{\partial H_1}{\partial \bar{\eta}_{k}}\big(G(z), \bar{G}(0)\big)\frac{\partial \bar{g}_k}{\partial \bar w}(0)}{H_1\big(G(z), \bar{G}(0)\big)} \nonumber\\
&-&\frac{\sum_{k=1}^{n}\frac{\partial H_2}{\partial \bar{\eta}_k}\big(G(z), \bar{G}(0)\big)\frac{\partial \bar{g}_k}{\partial \bar w}(0)}{H_2\big(G(z), \bar{G}(0)\big)}+\overline\varphi'(0).
\end{eqnarray}
Note that $H_1\big(G(z), \bar{G}(w)\big),  H_2\big(G(z), \bar{G}(w)\big)$  are two polynomials in $G(z)$ and $\bar{G}(w)$, $H_1\big(G(z), \bar{G}(w)\big)\neq 0$ everywhere on $U_1\times U_2$, and $H_2\big(G(z), \bar{G}(w)\big)\neq 0$ everywhere on $\mathbb{G}_{n}\times\mathbb{G}_{n}$. Then the right hand side of equation (\ref{3.7}) is a well-defined holomorphic rational function in $g_1, \ldots, g_n$.
Following the similar discussion, for any positive integer $\delta$, and for $z$ near 0, we get
\begin{equation}\label{3.8}
F(z) \cdot D^{\delta}(\bar{F}(0)) = Q_\delta(g_1, \ldots, g_n), ~~\delta=1, 2, \cdots,
\end{equation}
where $Q_\delta(g_1, \ldots, g_n)$ is a holomorphic rational function in $g_1,\ldots, g_n$.

Now denote $\mathcal{V}:=Span_{\mathbb{C}}\{D^{\delta}(\bar{F}(w))|_{w=0}\}_{\delta\geq1}$ be a vector subspace of $\mathbb{C}^m$. Since $F$ is nonconstant by our assumption, $\mathcal{V}$ cannot be a zero space. Now we let $\{D^{\delta_{j}}(\bar{F}(w))|_{w=0}\}_{j=1}^d$ be a basis for $\mathcal{V}$. Because $\bar{F}(w)$ is anti-holomorphic on $D$ and $\bar{F}(0)=0$, for any $w$ near 0, by the Taylor expansion we have:
\begin{equation*}
\bar{F}(w)=\sum_{\delta\geq1}\frac{D^{\delta}(\bar{F}(0))}{\delta!} {\bar w}^\delta \in \mathcal{V},
\end{equation*}
where $D^{\delta}$  is defined by (\ref{3.y}). Then for a small neighborhood $U_0$ of $0$, we have $\bar{F}(U_0)\subseteq\mathcal{V}$.

On the other hand, Taking the vectors $\{V_j\}_{j=1}^{m-d}$ as the basis of $\mathcal{V}^{\bot}$ on $\mathbb{C}^m$. Then we get
\begin{equation*}
F(z)\cdot V_j=0, ~~j=1,\ldots, m-d.
\end{equation*}
Combining with (\ref{3.8}), we can obtain a non-degenerate linear equations
\begin{equation*}
\left(f_1, \cdots, f_m
\right)
\left(
\begin{array}{c}
   D^{\delta_{1}}(\bar{F}(w))|_{w=0} \\
   \vdots\\
   D^{\delta_{d}}(\bar{F}(w))|_{w=0}\\
   V_1\\
   \vdots\\
   V_{m-d}
  \end{array}
          \right)^T
 = \left(
  \begin{array}{c}
  Q_{\delta_1}\\
   \vdots\\
  Q_{\delta_d}\\
  0\\
  \vdots\\
  0
 \end{array}
          \right).
\end{equation*}
It is obvious that each complement $f_j$ of $F$ can be linearly expressed by $\{Q_{\delta_j}\}_{j=1}^d$ by Gramer's rule, i.e., we can write
\begin{equation*}
f_j=\tilde{Q}_{j}(g_1, \ldots, g_n), ~~j=1, \cdots, m,
\end{equation*}
where $\{ \tilde{Q}_{j}(g_1, \ldots, g_n) \}_{j=1}^{m}$ are holomorphic rational functions in $g_1, \ldots, g_n$,
and then they are holomorphic Nash-algebraic functions in $g_1, \ldots, g_n$.

\textbf{Step 3.} Now we consider two cases to achieve the contradiction.

Let $\mathfrak{R}$ be the field of rational functions in $z$ over $D$. Consider the field extension
\begin{equation*}
\widetilde{\mathfrak{R}}=\mathfrak{R}(g_1, \ldots, g_n),
\end{equation*}
i.e., the smallest subfield of meromorphic function field over $D$ containing rational functions and $g_1, \ldots, g_n$. Denote  $trdge (\widetilde{\mathfrak{R}}/\mathfrak{R})$ be the transcendence degree of the field extension $\widetilde{\mathfrak{R}}/\mathfrak{R}$.

\textbf{Case 1.} If $trdge (\widetilde{\mathfrak{R}}/\mathfrak{R})=0$, i.e., $g_1, \ldots, g_n$ are all holomorphic Nash-algebraic functions in $z$.

From the above argument, we get that $f_j$ $(j=1, \ldots, m)$ are also Nash-algebraic functions in $z$.  Together with the following equation
\begin{equation*}
\begin{aligned}
\exp\big(\sum_{i=1}^{m}f_i(z)\bar{f_i}(w)\big)&=e^{\varphi(z)+\overline\varphi(w)}K_{\mathbb{G}_n}(G(z), {G}(w))\\
&=e^{-(\varphi(0)+\overline\varphi(0))}\frac{K_{\mathbb{G}_n}(G(z), {G}(w))}{K_{\mathbb{G}_n}(G(z), {G}(0))K_{\mathbb{G}_n}(G(0), {G}(w))}
\end{aligned}
\end{equation*}
and Lemma \ref{lem2.4}, we obtain that $F$ is a constant map.

\textbf{Case 2.} If $trdge (\widetilde{\mathfrak{R}}/\mathfrak{R}):=l>0$, i.e., $g_1, \ldots, g_n$ are not all holomorphic Nash-algebraic functions in $z$.

One can choose, without loss of generality, that $g_1, \ldots, g_l$ $(l\leq n)$ is the maximal algebraic independent subset in $\widetilde{\mathfrak{R}}$.
Then
\begin{equation*}
trdge \big(\widetilde{\mathfrak{R}}/\mathfrak{R}(g_1, \ldots, g_l)\big)=0.
\end{equation*}
Thus any element in $\{g_{l+1},\ldots, g_n\}$ is holomorphic Nash-algebraic function in $z, g_1, \ldots, g_l$. In other words, there exists a small neighborhood $V$ of $0$ such that for $\{g_i\}_{i=l+1}^{n}$,
we have some holomorphic Nash-algebraic functions $\{\hat{g}_{i}(z, X)\}_{i=l+1}^{n}$
in the neighborhood $\hat{V}$ of $\{(z, g_1,\ldots, g_l)|z \in V\}\subseteq \mathbb{C} \times \mathbb{C}^{l}$ such that
\begin{eqnarray*}
&&g_i(z)=\hat{g}_{i}(z, g_1, \ldots, g_l), ~~i=l+1, \ldots, n,
\end{eqnarray*}
for all $z\in V$, where $X=(X_1, \ldots, X_l)$. Then by the step two, there exist holomorphic Nash-algebraic functions $\{\hat{f_i}(z, X)\}_{i=1}^m$ on $\hat{V}$ such that
\begin{eqnarray*}
&&f_i(z)=\hat{f_i}(z, g_1, \ldots, g_l), ~~i=1, \ldots, m,
\end{eqnarray*}
for all $z\in V$.

Denote $\hat{G}(z, X)= (\hat{g}_{{l+1}}(z,X), \ldots, \hat{g}_{n}(z, X))$. Then, by $(\ref{3.2})$ and $(\ref{har})$, we define a function on $\hat{V}\times V$ as follows:
\begin{equation*}
\begin{aligned}
\psi(z, X, w)=&\sum_{i=1}^{m}\hat{f_i}(z, X)\bar{f_i}(w)-\log K_{\mathbb{G}_n}\big((X, \hat{G}(z, X)), {G}(w)\big)\\
&+\log K_{\mathbb{G}_n}\big((X, \hat{G}(z, X)), {G}(0)\big)+\overline\varphi(0)-\overline\varphi(w).
\end{aligned}
\end{equation*}
Then $\psi(z, g_1, \ldots, g_l, w)\equiv0$ on $V$. Now we claim that $\psi(z, X, w)\equiv0$ on $\hat{V}\times V$.

Define
\begin{equation*}
\ \phi(z, X, w)=\frac{\partial \psi}{\partial \bar w}(z, X, w).
\end{equation*}
We need only to prove that $\phi(z, X, w)\equiv 0$ on $\hat{V}\times V$.

Otherwise, then there exists a neighborhood $V_0$ of $0\in V$ such that $\phi(z, X, w_0)\not\equiv 0$. For fixed $w_0\in V_0$, $\phi(z, X, w_0)$ is holomorphic Nash-algebraic function in $(z, X)$.
Assume that its annihilating function is
\begin{equation*}
P(z, X, t)=a_d(z, X)t^d+\cdots+a_0(z, X),
\end{equation*}
where $a_0(z, X)\not\equiv0$ on $\hat{V}$, and $\{a_i(z, X)\}_{i=0}^d$ are holomorphic polynominals in $(z, X)$. Note that $\psi(z, g_1, \ldots, g_l, w_0) = 0$ on $V$. Then $ \phi(z, g_1, \ldots, g_l, w_0)= 0$ on $V$. Hence,
\begin{equation*}
P(z, g_1, \ldots, g_l, \phi(z, g_1, \ldots, g_l, w_0))=P(z, g_1, \ldots, g_l, 0)=a_0(z, g_1, \ldots, g_l)=0.
\end{equation*}
That is, we get $\{g_1, \ldots, g_l\}$ are algebraic dependent over $\mathfrak{R}$, which is a contradiction.
Therefore, $\psi(z, X, w)\equiv0$ for all $(z, X)\in \hat{V}$, $w \in V_0$, and then we have $\psi(z, X, w)\equiv0$ for $(z, X, w)\in \hat{V}\times V$.

Now we have the following equality
\begin{equation*}
\begin{aligned}
\sum_{i=1}^{m}\hat{f_i}(z, X)\bar{f_i}(w)=&\log K_{\mathbb{G}_n}\big((X, \hat{G}(z, X)), {G}(w)\big)\\
&-\log K_{\mathbb{G}_n}\big((X, \hat{G}(z, X)), {G}(0)\big)-\overline\varphi(0)+\overline\varphi(w),
\end{aligned}
\end{equation*}
where $(z, X, w)\in \hat{V}\times V$.
In particular, for some fixed $z$ and $w$, the left hand side of this equation is not identically equal zero on $\hat{V}\times V$. Indeed, if
\begin{equation*}
\sum_{i=1}^{m}\hat{f_i}(z, X)\bar{f_i}(w)\equiv 0,
\end{equation*}
by setting $w={z}$, we have
\begin{equation*}
\sum_{i=1}^{m}|f_i(z)|^2=\sum_{i=1}^{m}\hat{f_i}(z, g_1, \ldots, g_l)\bar{f_i}({z})\equiv 0.
\end{equation*}
This implies that $\{f_i(z)\}_{i=1}^m$ are constant maps, which contradicts with the previous assumption.

Next consider the following equation
\begin{equation}\label{3.9}
\exp\left(\sum_{i=1}^{m}\hat{f_i}(z, X)\bar{f_i}(w)\right)=e^{-(\varphi(0)+\overline\varphi(0))}\frac{K_{\mathbb{G}_n}\big((X, \hat{G}(z, X)), {G}(w)\big)}{K_{\mathbb{G}_n}\big((X, \hat{G}(z, X)), {G}(0)\big)K_{\mathbb{G}_n}\big(G(0), {G}(w)\big)}.
\end{equation}
Note that $\sum_{i=1}^{m}\hat{f_i}(z, X)\bar{f_i}(w)$ is a nonconstant holomorphic Nash-algebraic function in $X$ for some fixed $z$ and $w$, and the right hand side is also a holomorphic Nash-algebraic function in $X$, which is a contradiction by Lemma \ref{lem2.4}.

Therefore $F$ must be constant. This completes the proof of Theorem \ref{Th1.1}.

\vskip 10pt

\noindent\textbf{Acknowledgments}\quad The project is supported by the National Natural Science Foundation of China (No. 11671306).

\addcontentsline{toc}{section}{References}
\phantomsection
\renewcommand\refname{References}

\clearpage
\end{document}